\documentclass[reqno]{amsart}
\usepackage{amsmath, amssymb, amsthm, epsfig}
\usepackage[hidelinks]{hyperref}
\usepackage{latexsym}
\usepackage{url}
\usepackage[mathscr]{euscript}

\usepackage{color}
\usepackage{fullpage}
\usepackage{setspace}

\def\today{\ifcase\month\or
	January\or February\or March\or April\or May\or June\or
	July\or August\or September\or October\or November\or December\fi
	\space\number\day, \number\year}

\newtheorem{problem}{Extremal Problem}

\newtheorem{theorem}{Theorem}

\newtheorem{lemma}[theorem]{Lemma}

\theoremstyle{definition}

\theoremstyle{remark}

\title{Fourier optimization and consequences of the generalized Riemann hypothesis}

\author[Quesada-Herrera]{Emily Quesada-Herrera}

\address{University of Lethbridge, Department of Mathematics and Computer Science, 4401 University Dr W, Lethbridge, AB T1K 3M4}
\email{emily.quesadaherrera@uleth.ca}

\date{\today}

\begin{document}

\allowdisplaybreaks
\numberwithin{equation}{section}
\begin{abstract}
    We give an exposition of some connections between Fourier optimization problems and problems in number theory. In particular, we present some recent conditional bounds under the generalized Riemann hypothesis, achieved via a Fourier optimization framework, on bounding the maximum possible gap between consecutive prime numbers represented by a given quadratic form; and on bounding the least quadratic non-residue modulo a prime number. This is based on joint works with Emanuel Carneiro, Andr\'es Chirre, Micah Milinovich, and Antonio Pedro Ramos.
\end{abstract}
\maketitle

\section{Introduction}
This article is an exposition of connections between Fourier analysis and number theory, and is partly based on the article \cite{CQ-H} (joint with Andr\'es Chirre) and the preprint \cite{CMQR} (joint with Emanuel Carneiro, Micah B. Milinovich, and Antonio Pedro Ramos). We present two results, conditional on the generalized Riemann hypothesis, to illustrate this connection. Our results are inequalities regarding problems in number theory, and Fourier analysis was used to optimize the value of the explicit constants that appear in these inequalities.

To state the problems under consideration and our results, consider first quadratic and cubic residues modulo a prime number $p$. An integer $a$ is a quadratic residue mod $p$ if there exists an integer $u$ such  that $u^2\equiv a$ mod $p$. If this is not the case, $a$ is a quadratic non-residue. Similarly, $a$ is a cubic residue if there exists an integer $u$ such that $u^3\equiv a$ mod $p$.

Our first problem involves the study of integers and prime numbers represented by quadratic forms. This is a classical problem in number theory, going back at least to Fermat - he realized, for instance, that a prime number $p$ is a sum of two squares if and only if $p\equiv 1$ mod $4$. A more subtle example is the quadratic form $u^2+27v^2$. A conjecture of Euler, proven by Gauss, says that a prime $p$ has the form $p=u^2+27v^2$ if and only if both $p\equiv 1$ (mod 3) and $2$ is a cubic residue mod $p$. The following is an illustration of one of our results in \cite{CQ-H}.
\begin{theorem}\label{thm:forms1}
    Assume the generalized Riemann hypothesis. Then, for all sufficiently large $x$, there is always a prime number $p$ of the form $p=u^2+27v^2$ between $x$ and $x+\frac{46}{25}\sqrt{x}\log x$.
\end{theorem}
Our second result is about the maximum possible size of the least quadratic non-residue modulo a prime number. In \cite{CMQR}, we show the following.
\begin{theorem}\label{thm:lqnr}
    Assume the generalized Riemann hypothesis. Let $n_p$ be the least quadratic non-residue modulo a prime $p$. If $p$ is sufficiently large, then
    \begin{equation*}
        n_p< \frac{61}{80} \log^2 p.
    \end{equation*}
\end{theorem}
Below we give more information on these problems and their background.

\section{Fourier transform}\label{sec:fourier}
The Fourier transform is certainly one of the most fundamental objects in mathematics and applied mathematics, as it is used to model a variety of oscillatory phenomena.
For an integrable function $f:\mathbb{R}\to \mathbb{C}$, we may define the Fourier transform of $f$ by
\begin{equation*}
    \widehat{f}(\xi):=\int_{\mathbb{R}} f(x)e^{-2\pi i \xi x }dx \ \  \ \ (\xi\in\mathbb{R}).
\end{equation*}
%
A Fourier optimization problem is one where we impose conditions on a function and its Fourier transform, and seek to optimize a given quantity. These problems often have surprising applications, such as in the theory of the Riemann zeta-function and the distribution of its zeros \cite{CCM, CChi, CChiM,CMR,CS,CPL, Q-H}, in bounding prime gaps \cite{CMS,CPL, CQ-H}, and in the theory of sphere packings \cite{CE,CKMRV, Vi}. A classical reference on using extremal functions in Fourier analysis is Vaaler's 1985 survey \cite{Vsurvey}.

\section{Riemann zeta-function and Dirichlet $L$-functions}\label{sec:riemann}
Understanding the behavior of the Riemann zeta-function and the distribution of its zeros is a crucial problem in number theory, being related to the distribution of primes. The Riemann zeta-function is defined in the half-plane $\text{Re}s>1$ by
\begin{equation}\label{def:zeta}
    \zeta(s)=\sum_{n=1}^\infty \frac{1}{n^s}=\prod_p \left(1-\frac{1}{p^s}\right)^{-1}.
\end{equation}
The \textit{Euler product} - the product over primes in \eqref{def:zeta} - is the starting point of the connection between this function and the distribution of primes.
It is possible to show that $\zeta(s)$ extends to a meromorphic function in the entire complex plane, with a simple pole at $s=1$, and that $\zeta(-2n)=0$ for all positive integers $n$. These are called the \textit{trivial zeros.} Furthermore, one can show that all non-trivial zeros of $\zeta(s)$ must lie on the \textit{critical strip} $\{s\in\mathbb{C}:0<\text{Re}s <1\}$. The Riemann hypothesis (RH), conjectured in 1859 \cite{RH}, states that all non-trivial zeros actually lie on the \textit{critical line}, $\text{Re}s = \tfrac{1}{2}$.
For further background on $\zeta(s)$ and primes, see \cite{DP,IK, MV2, Tit0}.

Dirichlet $L$-functions are one generalization of $\zeta(s)$, and were originally introduced to study primes in arithmetic progressions. Instead of giving a full definition, let us see a concrete example. Consider the Legendre symbol, defined by
\begin{equation*}
    \chi(a)=\left(\frac{a}{p}\right):=\left\{
    \begin{array}{ll}
        0 & \text{ if } p\mid a; \\
        1 & \text{ if } p\nmid a \text{ and } a \text{ is a quadratic residue mod } p;\\
        -1 & \text{ if } p\nmid a \text{ and } a \text{ is a quadratic non-residue mod } p.
    \end{array}
    \right.
\end{equation*}
This is an example of a \textit{Dirichlet character} modulo $p$, and in this case, it codifies information about quadratic non-residues. We consider the associated Dirichlet $L$-function, defined initially (and then analytically continued), for $\text{Re}s>1$, by
\begin{equation}\label{eq:dirichlet}
    L(s,\chi)= \sum_{n=1}^\infty \frac{\chi(n)}{n^s}=\prod_p \left(1-\frac{\chi(p)}{p^s}\right)^{-1}.
\end{equation}
The \textit{generalized Riemann hypothesis} (GRH) conjectures that all non-trivial zeros (i.e, those in the critical strip) of $L(s,\chi)$ also have real part $\tfrac{1}{2}$. See, for instance, \cite{DP, IK} for details and a full definition.
\subsection{A link between worlds}
A bridge that connects Fourier analysis with number theory is the Guinand-Weil explicit formula, and it is an ingredient in many of the aforementioned works. It connects an arbitrary function $h$, its Fourier transform $\widehat{h}$, prime numbers, and zeros of the Riemann zeta-function or an $L$-function:
\begin{lemma}[Guinand-Weil explicit formula] \label{GW_lemma} Let $h(s)$ be analytic in the strip $|\mathrm{Im} \, s| \le \tfrac12+\varepsilon$ for some $\varepsilon>0$, and assume that $|h(s)| \le C(1+|s|)^{-(1+\delta)}$ for some $\delta>0$ and $C>0$. Let $\chi$ be a primitive Dirichlet character modulo $q$. Then
\begin{eqnarray*}\label{GW}
\sum_{\rho_{\chi}} h\!\left(\frac{\rho_{\chi}-\tfrac12}{i}\right) &= \widehat{h}(0) \dfrac{\log (q/\pi)}{2\pi}   +  \dfrac{1}{2\pi}  \displaystyle\int_{-\infty}^{\infty} h(u) \, \mathrm{Re} \,  \frac{\Gamma'}{\Gamma}\left(\frac{2 - \chi(-1)}{4} + \frac{iu}{2}\right)  du
\\
& - \dfrac{1}{2\pi}  \displaystyle\sum_{n\geq2}\frac{\Lambda(n)}{\sqrt{n}} \left\{  \chi(n) \, \widehat{h}\!\left( \frac{\log n}{2\pi} \right) + \overline{\chi(n)}\,\widehat{h}\!\left( -\frac{\log n}{2\pi} \right)  \right\},
\end{eqnarray*}
where the sum on the left-hand side runs over the non-trivial zeros $\rho_{\chi}$ of $L(s,\chi)$, and $\Lambda(n)$ is the von Mangoldt function defined to be $\log p$ if $n=p^k$, $p$ a prime and $k\ge 1$, and zero otherwise.
\end{lemma}
Here, the information about prime numbers is codified in the Von Mangoldt function $\Lambda(n)$, and, in our case above, $\chi(n)$ contains information about quadratic residues.
While this gives a bridge between the fields, formulating a suitable Fourier optimization problem is a subtle and delicate matter.
\section{Prime gaps and quadratic forms}\label{sec:forms}
Due to the connection between primes numbers and zeros of the Riemann zeta-function and other $L$-functions, it is useful to study the consequences of RH and GRH, as it may allow us to better understand the limits of our current techniques when tackling a certain problem. Let $p_n$ be the $n$-th prime number. Assuming RH, Cram\'er showed that the gap between consecutive prime numbers satisfies
\begin{equation}\label{eq:cramer}
    p_{n+1}-p_n\le (C+o(1))\sqrt{p_n}\log p_n
\end{equation}
for an absolute constant $C>0.$ The order of magnitude here has not been improved in the past 90 years, and the efforts have been in reducing the explicit value of $C$. Carneiro, Milinovich and Soundararajan \cite{CMS} developed a Fourier optimization strategy to tackle this. They introduce the following Fourier optimization problem.
\begin{problem}
 Define $\mathcal{A^+}$ to be the class of even and continuous functions $F:\mathbb{R}\rightarrow \mathbb{R}$, with $F\in L^1(\mathbb{R})\setminus\{0\}$. For a parameter $1\le A < \infty,$ find
		\begin{equation}\label{eq:cms}
		\mathcal{C^+}(A):= \sup_{F\in \mathcal{A^+}} \frac{F(0) -A\int_{[-1,1]^c}(\widehat{F}(t))_+ \,d t}{\|F\|_1}.
		\end{equation}
\end{problem}
Here, we denote the $L^1$ norm $\|F\|_1:=\int_\mathbb{R} |F(x)| \, dx$, and $\widehat{F}(t)_+=\max\{\widehat{F}(t),0\}$. They show that, assuming RH,
we may take the constant $C=\mathcal{C^+}\left(36/11\right)^{-1}$ in the bound \eqref{eq:cramer},
 giving this connection between a problem in Fourier analysis and the problem of bounding prime gaps. The specific parameter $A=\frac{36}{11}$ arises from other number-theoretic considerations. This Fourier optimization problem can be studied on its own, and they show several qualitative and quantitative properties. For example, they show that $1\le \mathcal{C^+}(A)\le 2$, for all $A\ge 1,$ and furthermore, for any $A>1$, there exists an extremal function $F$ that yields the maximum value in \eqref{eq:cms}.

 While this problem can be stated in accessible terms, actually constructing an extremal function $F$, or even finding the exact value of $\mathcal{C^+}(A)$ for a given parameter $A$, is a difficult problem. To make the constant in \eqref{eq:cramer} as small as possible, the question is then to give good lower bounds for $\mathcal{C^+}(36/11)$, by means of constructing a good explicit test function $F$. They find a specific example that gives the bound $\mathcal{C^+}(A)> \frac{25}{21}$, improving the then best-known admissible value of $C$ in \eqref{eq:cramer} to $p_{n+1}-p_n\le \left(\frac{21}{25}+o(1)\right)\sqrt{p_n}\log p_n$, and this is still the best approach to date to attack this problem. The numerics where slightly refined by Chirre, Pereira and De Laat \cite{CPL} by using strong computational techniques to construct test functions, via semidefinite programming. They replaced $\frac{21}{25}$ by $0.8358$, and also used the method of Carneiro, Milinovich and Soundararajan to study gaps between primes in an arithmetic progression. Returning to \cite{CMS}, they also showed the following explicit version, with a similar but slightly larger constant: for all $x>2$,  there is always a prime number between $x$ and $x+\frac{22}{25}\sqrt{x}\log x$.

In the article \cite{CQ-H}, we combine tools from Fourier analysis, analytic number theory and algebraic number theory to obtain several new estimates regarding integers and primes represented by quadratic forms. We have two main themes that are ubiquitous in this investigation. First, we use the well-known theme that propositions about quadratic forms can be stated in two other equivalent languages: ideals of number fields and lattices. We use all three points of view to our advantage in different parts of the article. Our second theme is the use of Fourier analysis, in the following way: we begin by finding a summation formula that connects our object of study with an arbitrary function and its Fourier transform; then, we choose an appropriate test function that recovers the desired information in an optimized manner. The use of a version of the Guinand-Weil explicit formula is one example of this (which is the focus of the present expository article), but not the only one. One of our results is an extension of the method of Carneiro, Milinovich, and Soundararajan to study gaps between primes represented by a given quadratic form. Consider, for instance, the quadratic form $q(u,v)=u^2+27v^2$, and let $p_{n,q}$ be the $n$-th prime number of the form $p=u^2+27v^2$.
The following is a corollary of our results:
\begin{theorem}
    Assume GRH. Then,
    \begin{equation*}
        p_{n+1,q}-p_{n,q} \le
        \left(\frac{6}{\mathcal{C}^+(28)} +o(1)
        \right) \sqrt{p_{n,q}}\log p_{n,q}.
    \end{equation*}
    Moreover, we have the lower bound $\mathcal{C}^+(28)\ge  1.0889.$
\end{theorem}
This immediately implies Theorem \ref{thm:forms1}. To achieve this, we work in the language of ideals of imaginary quadratic fields, where we have the machinery of Hecke $L$-functions- another generalization of Dirichlet $L$-functions to this language- for which we can use an explicit formula similar to Lemma \ref{GW_lemma}. To obtain the explicit bound for $\mathcal{C}^+(28)$, we use the following test function in the optimization problem \eqref{eq:cms}:
\begin{equation*} \label{eq:finalInequalityF3}
	 F(x)=H\left(\frac{x}{0.98644}\right),
\end{equation*}
	where
 \begin{equation*}
     H(x)= \cos(2\pi x)\left(
	\frac{68}{1-16x^2} + \frac{5}{9-16x^2}+\frac{1}{25-16x^2}
	\right).
 \end{equation*}
 In Section \ref{sec:numerical}, we explain how to find such test functions numerically.
%

\section{Least quadratic non-residue}\label{sec:lqnr}
Let us now turn our attention to the problem of bounding $n_p$, the least quadratic non-residue modulo a prime number $p$. For example, $1=1^2$ is always a quadratic residue. How about 2? Who is the smallest non-residue?
Ankeny showed in 1952 that, assuming GRH, we have the bound $n_p\le (C+o(1))\log^2 p$, so that we can always find a non-residue that is much smaller than $p$.
For further background, see \cite{CMQR} and the references therein. In \cite{CMQR}, we consider the following optimization problem:
\begin{problem}
Find
\begin{equation}\label{eq:prob2}
    \mathcal{C}(1):= \sup_{0\neq F\in\mathcal{A}} \frac{2\pi}{\|F\|_1}
    \left(
    \int_{-\infty}^0 \widehat{F}(t) e^{\pi t}\, d t -
    \int_{0}^\infty |\widehat{F}(t)| e^{\pi t}\, d t
    \right),
\end{equation}
where the supremum is taken over the class of functions $\mathcal{A}=\{F:\mathbb{R}\to\mathbb{C};\, F\in L^1(\mathbb{R}):\, \widehat{F}$ is real-valued $
\}$.
\end{problem}
We show the following connection:
\begin{theorem}
    Assume GRH. Let $n_p$ be the least quadratic non-residue modulo a prime $p$. Then
    \begin{equation*}
        n_p\le (\mathcal{C}(1)^{-2} + o(1)) \log^2 p.
    \end{equation*}
    Furthermore, we have the bounds $1.146 < \mathcal{C}(1) < 1.148.$
\end{theorem}
This immediately implies Theorem \ref{thm:lqnr}, and gives an improvement on the previous best asymptotic upper bounds under GRH, given by Lamzouri, Li, and Soundararajan \cite{LLS}. An interesting feature of this result is that, since we have both lower bounds (obtained by constructing good functions $f$) and upper bounds (which are more subtle and require finer theoretical analysis of the optimization problem), we can see that we have essentially arrived at the limit of this method.
\section{Numerically optimizing the bounds}\label{sec:numerical}
Once we arrive at a Fourier optimization problem, it can be difficult to find an exact solution. The next step is to find good bounds to approximate it- that is, find a test function $F$ that does a good job in the search for the supremum.
A way to approach this is to start with a good family of test functions that parametrize our search space, and then use numerical optimization methods to find a good candidate within that family. One possible family is functions of the form $P(x)e^{-\pi x^2}$, where $P$ is a polynomial: this family was used in \cite{CPL}, where candidates where found via semidefinite programming (see \cite{semi} for background on semidefinite programming). In other cases, bandlimited functions- those whose Fourier transform have compact support- are useful, such as in \cite{CQ-H}.

\begin{figure}[t]
        \centering
        \includegraphics[width=5cm]{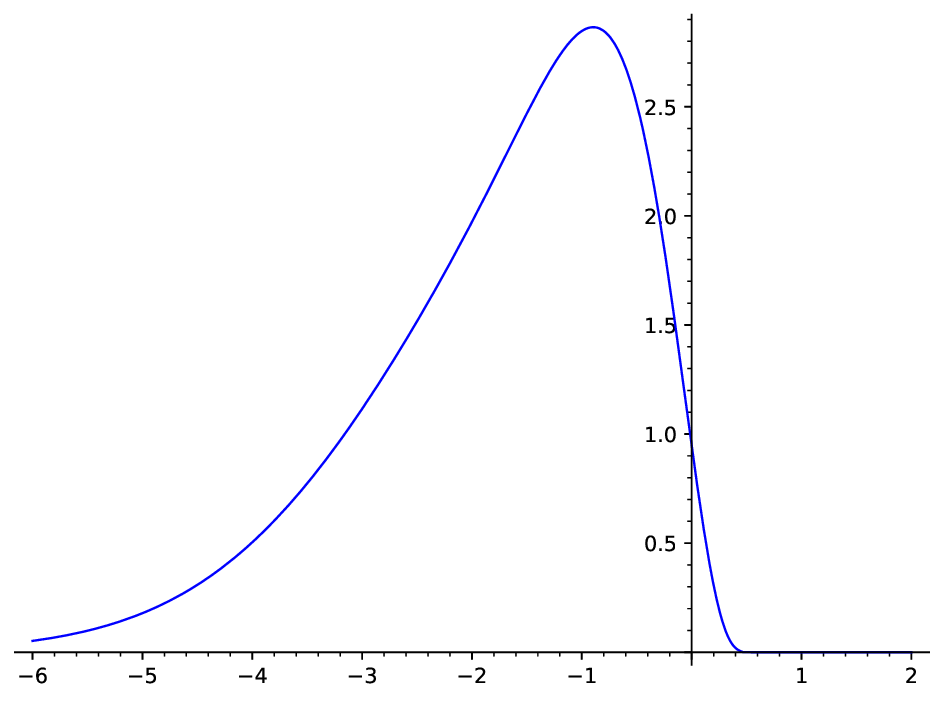}
        \caption{Plot of $\widehat{F}(t)/\widehat{F}(0)$ defined in \eqref{eq:F}, a good test function for the optimization problem \eqref{eq:prob2}.}
        \label{fig:my_label}
    \end{figure}

To illustrate this process, let us consider the extremal problem \eqref{eq:prob2}. To maximize the numerator, we should try to concentrate most of the mass of $\widehat{F}(t)$ on the negative numbers. With this in mind, we consider the Fourier transform pairs
\begin{equation}\label{20240411_14:28}
F_n(x)= \frac{(2n-1)!}{\pi(1 + 2ix)^{2n}}; \qquad \widehat{F_n}(t) = -(\pi t)^{2n-1} \,e^{\pi t}\,{\bf 1}_{\mathbb{R}_-}(t),
\end{equation}
and note that $F_n(x)\in \mathcal{A}$ for all positive integers $n$. Here, ${\bf 1}_{\mathbb{R}_-}(t)$ is the indicator function of the negative numbers, defined as $1$ if $t<0$ and 0 otherwise.
We now consider the following family, composed by dilations and translations of linear combinations of functions as in \eqref{20240411_14:28}, defined by:

\begin{equation}\label{eq:lowers}
 \widehat{F}(t)= g\left(\frac{\pi t - c}{a}
 \right)
 ,\text{ where } \ g(t)=\sum_{n=1}^N b_n\, t^{2n-1} \, e^{t}\, {\bf 1}_{\mathbb{R}_-}(t).
\end{equation}
Here, $b_n, c\in \mathbb{R}$ and $a>0$. Now we use numerical optimization methods to find good parameters. With the principal axis method of Brent \cite{Br},
%
we find the function
\begin{equation}\label{eq:F}
    \widehat{F}(t)= g\left(\frac{t - 0.47}{0.42}\right),
\end{equation}
where
\begin{equation*}
        g(x)= 2e^x  \left(0.0006 x^7+0.0005 x^5+ x^3+0.0405 x\right) {\bf 1}_{\mathbb{R}_-}(x).
\end{equation*}
See Figure \ref{fig:my_label} for a plot of $\widehat{F}(t)$.
This yields the bound $1.143 <\mathcal{C}(1)$, which is close to what we claim. In \cite{CMQR}, we use similar, more complicated examples to get our result.
\section*{Acknowledgements}
The author was supported by a PIMS-Simons postdoctoral fellowship at the University of Lethbridge. The author thanks the anonymous referee for helpful comments.

\bibliographystyle{abbrv}
\bibliography{fogrh}

\end{document}